\documentclass{amsart}	
\usepackage{amsmath,amssymb,amsthm,amscd}
\usepackage{graphicx,subfigure}
\usepackage{color}
\usepackage{pb-diagram} 
\usepackage{tabularx}
\usepackage{multirow}
\usepackage{longtable}
\usepackage{hyperref}
\usepackage[foot]{amsaddr}
\usepackage{enumitem}
\usepackage{comment}
\usepackage{phonetic}

\usepackage[utf8]{inputenc}

\usepackage{geometry}
\newcommand{\blue}[1]{\bgroup\color{blue}{#1}\egroup}
\newcommand{\red}[1]{\bgroup \egroup}
\newcommand{\magenta}[1]{\bgroup\color{magenta}{#1}\egroup}
\newcommand{\green}[1]{\bgroup\color{green}{#1}\egroup}

\newtheorem{theorem}{Theorem}[section]

\newtheorem{lemma}[theorem]{Lemma}

\theoremstyle{remark} 
\newtheorem{remark}[theorem]{Remark}

\numberwithin{equation}{section}

\renewcommand\epsilon{\varepsilon}

\newcommand{\cD}{{\mathcal D}}

\newcommand{\RR}{{\mathbb R}}
\newcommand{\CC}{{\mathbb C}}

\newcommand{\TT}{{\mathbb T}}
\newcommand{\ZZ}{{\mathbb Z}}

\newcommand{\abs}[1]{|{#1}|}

\newcommand{\norm}[1]{\left\|{#1}\right\|}

\newcommand{\aver}[1]{\langle{#1}\rangle}

\newcommand\wnorm[3]{{\ \max\left\{\norm{\etaL}_{#2}, \frac{1}{\gamma #3^\tau}  \norm{\etaN}_{#2}\right\}}}
\newcommand\wnormeta{{\wnorm{\eta}{\rho}{\delta}}}

\newcommand\qwnormeta{{{\wnormeta}^2}}

\newcommand{\Lie}[1]{\mathfrak{L}_{#1}}

\newcommand\Xh{{X_h}}

\newcommand{\Dioph}[3]{\cD_{#1,#2}^{#3}}

\newcommand\Loper{{\mathfrak{L}_{\omega}}}
\newcommand\Roper{{\mathfrak{R}_{\omega}}}
 
\newcommand{\comp}{{\!\:\circ\!\,}}

\newcommand\K{{K}}

\newcommand\Kn{{\bar{\bar K}}}

\newcommand\En{{\bar{\bar E}}}

\newcommand\ttop{{\!\top}}

\newcommand\DK{{\Dif\K}}

\newcommand\DKn{{\Dif\Kn}}
\newcommand\DKT{({\Dif\K})^{\ttop}}

\newcommand\DKnT{({\Dif\Kn})^{\ttop}}

\renewcommand\L{{\DK}}
\newcommand\LT{{\DKT}}

\newcommand\Ln{{\Dif\Kn}}

\newcommand\N{{N}}
\newcommand\NT{{\N}^\ttop}

\newcommand\Nn{{\bar {\bar{\N}}}}
\newcommand\NnT{{\bar{\bar\N}}\phantom)\!\!^{\ttop}}

\newcommand\B{{B}}

\newcommand\Bn{{\bar{\bar B}}}

\newcommand\T{{T}}

\newcommand\Tn{{\bar{\bar T}}}

\newcommand\etaL{\eta^{{\rm D}\!K}}
\newcommand\tetaL{{\tilde \eta}^{{\rm D}\!K}}
\newcommand\etaN{\eta^N}

\newcommand\xiL{\xi^{{\rm D}\!\K}}
\newcommand\xiN{\xi^N}

\newcommand\bL{b^{{\rm D}\!K}}
\newcommand\tbL{{\tilde b}^{{\rm D}\!K}}
\newcommand\bN{b^N}

\renewcommand\to{\rightarrow}

\newcommand\ii{\mathrm{i}}

\newcommand\Dif{ {\mbox{\rm D}} }

\newcommand\cmani{{\mathcal{U}}}

\renewcommand\H{h}                   % Hamiltonian
\newcommand\DeltaK{\Delta K}

\newcommand\cteXH{c_{\text{\tiny $\Xh$}}}
\newcommand\cteDXH{c_{\text{\tiny ${\rm D} \Xh$}}}
\newcommand\cteDDXH{c_{\text{\tiny ${\rm D}^2 \Xh$}}}
\newcommand\cteDXHT{c_{\text{\tiny $({\rm D} \Xh)^\ttop$}}}

\newcommand\sigmaDK{\sigma_{\mbox{{\tiny ${\rm D} K$}}}}
\newcommand\sigmaDKT{\sigma_{\mbox{{\tiny $({\rm D} K)\!^\ttop$}}}}

\newcommand\sigmaN{\sigma_{\mbox{{\tiny $\N$}}}}
\newcommand\sigmaNT{\sigma_{\mbox{{\tiny $\NT$}}}}

\newcommand\sigmaB{\sigma_B}
\newcommand\sigmaTinv{\sigma_{\mbox{{\tiny {$\aver{T}^{\mbox{-}1}$}}}}}

\newcommand\Csym{C_{\mathrm{sym}}}

\newcommand\CtheoE{\mathfrak{C}}
\newcommand\CtheoDeltaK{\mathfrak{C}_{\DeltaK}}

\newcommand\CxiL{C_{\xi^{\text{\tiny $L$}}}}
\newcommand\CDeltaKn{C_{\Delta \Kn}}

\newcommand\CDeltaBn{C_{\Delta \Bn}}
\newcommand\CDeltaNn{C_{\Delta \Nn}}
\newcommand\CDeltaNnT{C_{\Delta \NnT}}
\newcommand\CDeltaiTn{C_{\Delta \aver{\Tn}^{\text{-}1}}}
\newcommand\etan{{\bar {\bar \eta}}}
\newcommand\etanL{{{\etan}^L}}
\newcommand\etanN{{\etan^N}}
\newcommand\QCetanN{Q_{\etanN}}
\newcommand\QCetanL{Q_{\etanL}}
\newcommand\CDeltaDKn{C_{\Delta {\rm D}\Kn}}
\newcommand\CDeltaDKnT{C_{\Delta({\rm D}\Kn)^\ttop}}

\newcommand\hCDelta{\hat{C}_\Delta}

\newcommand\Hcar{H_{\rm C}}
\newcommand\Hdel{H_{\rm D}}
\newcommand\Hred{H}
\newcommand\Del{D}

\newcommand\Np{{n}}

\allowdisplaybreaks

\begin{document}
\title[]{
Sun-Jupiter-Saturn system may exist: a  
verified computation of 
quasiperiodic solutions for the planar three body problem
}
\date{\today}

\author{Jordi-Llu\'is Figueras$^{\mbox{1}}$}
\address[1]{Department of Mathematics, Uppsala University, 
Box 480, 751 06 Uppsala, Sweden}
\email{figueras@math.uu.se}

\author{Alex Haro$^{\mbox{2,3}}$}
\address[2]{Departament de Matem\`atiques i Inform\`atica, Universitat de Barcelona,
Gran Via 585, 08007 Barcelona, Spain.}
\address[3]{Centre de Recerca Matem\`atica, Edifici C, Campus Bellaterra, 08193 Bellaterra, Spain}
\email{alex@maia.ub.es}

\begin{abstract}
{In this paper, we present evidence of the stability of a simplified model of the Solar System, a flat (Newtonian) Sun-Jupiter-Saturn system with realistic data: masses of the Sun and the planets, their semi-axes, eccentricities and (apsidal) precessions of the planets close to the real ones. The evidence is based on convincing numerics that a KAM theorem can be applied to the Hamiltonian equations of the model to produce quasiperiodic motion (on an invariant torus)  with the appropriate frequencies.
To do so, we first use KAM numerical schemes to compute translated tori to continue from the Kepler approximation (two uncoupled two-body problems) up to the actual Hamiltonian of the system, for which the translated torus is an invariant torus. Second, we use KAM numerical schemes for invariant tori to refine the solution giving the desired torus. Lastly, the convergence of the KAM scheme for the invariant torus is (numerically) checked by applying several times a KAM iterative lemma, from which we obtain that the final torus (numerically) satisfies the existence conditions given by a KAM theorem.
} 
\end{abstract}

\maketitle

\section{Introduction}

In \cite{Newton} Newton
deduced the equations for the 
motion
of planets and solved the 2 body problem: bounded orbits
follow Kepler's motions (spin in ellipses with one focus on the center of mass), 
and unbounded ones are parabolae 
or hyperbolae. 
Then Newton (Book 3, Proposition XIII, Theorem XIII)
admits that observed planetary motion Jupiter does not fit the equations, and explains it by noticing 
that Saturn's influence can not be neglected. 
Since then, one of the most important problems in mathematics 
has been understanding the dynamics of the 3 (or higher) body problem. Many researchers have 
pursued this question and realized in different temporal 
stages that there are two (among others) important questions: the stability of the solutions - do 
planets orbit around the Sun in a quasiperiodic motion ad perpetuum? -; and the existence of chaos. 
This dichotomy was started by the pioneer work of Poincar\'e in \cite{Poincare_Book3BP}. In this 
paper we are interested in the 
stability problem.

Several steps forward in time and we encounter a fundamental advance towards solving the stability
problem. In 1954 Kolmogorov \cite{Kolmogorov54} presented a methodology for proving the existence 
of Lagrangian invariant tori in Hamiltonian systems of $n$ degrees of freedom close to integrable 
ones. Then Arnold \cite{Arnold63a, Arnold63b} and Moser \cite{Moser62} further explored this 
and the KAM theory was officially born. 
Since then a lot results have been produced, covering Lagrangian and lower dimensional tori, infinite 
dimensional systems, dissipative systems, etc. 
For the interested reader, we 
refer to the books \cite{BroerHS96, Llave01, Chierchia03}, 
and the popular book \cite{Dumas14}.

It was clear from the very beginning that the 3 body problem posed several obstacles that 
other Hamiltonians don't have. The integrable problem (Kepler's Hamiltonian)
doesn't have all the frequencies that the full problem has: in a general four degrees 
of freedom Hamiltonian invariant tori have four dimensions with four frequencies; while in 
Kepler's Hamiltonian they have two dimensions with two frequencies.
In KAM terminology it is said that 
the system is degenerate, and 
then
the full problem has several time-scale frequencies: 
the fast frequencies that correspond to the spinning of the planets around the Sun, 
and the slow frequencies that correspond to the spinning of their orbital ellipses (precession motion) 
and, in the spatial problem, 
the changes in the inclination of the rotation planes  (inclination motion).

A crucial advance was performed by Arnold in \cite{Arnold62, Arnold63b} 
where he proved the persistence of 
quasiperiodic motion for the planar three body problem for a ratio of the semi-major axis close to zero.  
The theory was later completed for the spatial $N$ body problem in remarkable works by 
Herman and  F\'ejoz \cite{Fejoz04}, and Chierchia and Pinzari \cite{ChierchiaP11}, among others.
In spite of the {\em fundamental} importance of 
all these theoretical results, they suffer
the {\em practical} inconvenience that the ratio of the semi-major axis or 
the size of masses of the planets (used as
the parameter measuring the distance to integrability) 
have to be ridiculously small. 
In fact, H\'enon \cite{Henon66} already 
took Arnold's paper and
checked that this size, in the simpler and non-degenerate {\em restricted} three body problem, 
is of order $10^{-333}$ (see the beautiful exposition of these facts in \cite{Laskar14}). 
After this result H\'enon asserted 
\footnote{From the English translation in \cite{Dumas14}}:
{\em ``Thus, these
theorems, although of a very great theoretical interest, do not seem applicable in their
present state to practical problems, where the perturbations are always much larger than the
thresholds [above].''}
This apparent lack of applicability of KAM theory to practical and physical problems led over time 
to some misunderstandings  (and laugther) about KAM theory 
but, as Dumas emphasizes in \cite{Dumas14}, H\'enon himself goes on to write: 
{\em 
``The numerical results we present here, and those obtained for other problems, indicate however that the 
[invariant] curves continue to exist for very strong perturbations, of 
the same order of magnitude as the leading term.''
}

This last observation by H\'enon is in fact what leads our research: the combination of qualitative
KAM results with computers. The idea is the use of the computers for getting initial data 
that can be then checked to fulfill the conditions of the taylored KAM theorems. Following this 
line of thought in combination with all the previous (classical) KAM methodology, based on performing 
canonical transformations on the Hamiltonian, 
there has been important advancement
towards the solution of the three body problem for realistic masses 
\cite{LaskarR95, Robutel95, LocatelliG05, LocatelliG07, 
CellettiC95}. 
More recently, in \cite{Castan2017thesis} a quantitative version of 
Arnold's KAM theorem have been applied to the plane
three-body problem to show, computer-assisted, the existence of quasiperiodic motion 
for a ratio of masses between the planets and the star 
that is close to $10^{-85}$ (this estimate accounts for a mass of the planets smaller 
than $10^{-24}$ times the mass of the electron).
In this paper, however, we propose to use another approach, based on the so-called parameterization method 
(see the seminal works \cite{Llave01,GonzalezJLV05}), looking directly for the 
parameterizations of the invariant tori, mitigating the curse 
of dimensionality.  In this approach, the KAM theorems are written in 
a posteriori format, so that the results are
suitable for numerical verification and, finally, for Computer-Assisted Proofs. This program 
has led to a rapid development of results \cite{HaroCFLM16, FiguerasHL17,HaroL19, FiguerasHL20, CallejaCLa}.

Applying KAM theory to the planar three body problem 
(for realistic parameters and ephemerides) is 
a very demanding problem, both mathematically and computationally. Hence, further steps must 
be performed to attack the problem. We have split this enterprise in three stages, of which the 
present paper is the centerpiece. 
Each stage deals with different questions and methods, so they could be 
of independent interest for different publics. Moreover, even though we have been thinking in 
the application to the three body problem, and specifically to the Sun-Jupiter-Saturn system, 
the pieces can be applied to other problems. The first stage, appearing in 
\cite{FiguerasHTheorem}, is a KAM theorem based on a (modified) parameterization method 
for Hamiltonian systems, with sharp 
control on the bounds and the Diophantine frequencies (with precedents in 
\cite{HaroL19,Villanueva17}). This first paper has two results that we use 
in this paper: the \emph{KAM Theorem} for verifying the existence of the invariant torus, 
and the \emph{Iterative Lemma} used for, giving an initial approximation with bounds on it, we perform several steps of the convergence scheme. This allows 
to refine the constants to be used later on the KAM Theorem.
The second stage, this paper, is a methodology to compute 
invariant tori in (close to degenerate) Hamiltonian systems with fast and slow 
time-scales, applied to numerically verify the existence of quasiperiodic 
solutions of the Sun-Jupiter-Saturn in the planar model with realistic masses and ephemerides. 
The last stage is \cite{FiguerasHValidation} 
where we present how the numerics from this paper and the KAM theorem from 
\cite{FiguerasHTheorem} are combined along with rigorous numerics for validating the results.

\subsection{Our results and their organization in the paper} 

We start  presenting the model we work on. 
It is a Hamiltonian with 3 degrees of freedom (the total angular momentum has been reduced) 
depending on a parameter $\mu$ that accounts for the 
masses, so that $\mu=  0$ corresponds to two uncoupled Kepler problems, and 
$\mu= \mu_0$ corresponds to the actual values of the masses of the planets (in our case $\mu_0= 10^{-3}$).
Then we discuss the numerical methods used. 
The goal is computing a 3 dimensional invariant torus for the observed values 
of the frequencies (for $\mu= \mu_0$).
Two of the frequencies are fast, and the other is slow and of the order of $\mu_0$. 
A fundamental obstacle we encounter is that the torus does not come by 
continuation from a 3 dimensional invariant torus for $\mu= 0$, because 
Kepler motions correspond to 2 dimensional tori, and the slow frequency
collapses to zero. Hence, we cannot apply a direct continuation technique of the {\em invariant tori} 
from the integrable problem since the problem is singular at $\mu= 0$.
Following the lines of thought of \cite{GonzalezHL13,GonzalezHL22}
we perform a continuation of {\em translated  tori},  
which are invariant tori for a modified Hamiltonian system to which we have added an extra term 
(a translation) that compensates the degeneracies of the actual problem. 
At $\mu= \mu_0$, the translation term should be zero.
At  this stage the torus is no longer degenerate, allowing the use of KAM numerical schemes on 
the actual problem for its
refinement. 
This leads to, after iterating several times 
the KAM numerical scheme, obtaining a very accurate approximation for the invariant torus.
Finally, with this approximation, we can run the Iterative Lemma
in \cite{FiguerasHTheorem}
several times and, lastly, the KAM Theorem so that 
all the bounds satisfy it and gives us the 
the existence of a nearby invariant torus, hence giving a numerical verification 
of the existence of quasiperiodic solutions close to the ephemerides of
Sun-Jupiter-Saturn configuration. 
Paraphrasing Hen\'on, the numerical results we present here indicate that the invariant 
torus exists for $\mu= \mu_0$.

\section{The planetary model and the problem} 

The planar $(1+\Np)$-body problem (the  Sun plus $\Np$ planets) in Poincar\'e heliocentric 
cartesian coordinates has Hamiltonian \cite{LaskarR95, ChierchiaP11} 
$\Hcar:\RR^{2n}\times\RR^{2n}\to \RR$ given by
\begin{equation}\label{eq: hamiltonian cartesian}
\begin{split}
\Hcar( x, y) & = \sum_{i=1}^\Np \left(\dfrac{\|y_i\|^2}{2m_i}-\dfrac{m_i}{\|x_i\|}\right)
+\mu \left(\sum_{i=1}^\Np\dfrac{\|y_i\|^2}{2} 
+ \sum_{1\leq i < j \leq \Np} \left(y_i\cdot y_j-\dfrac{m_im_j}{\|x_i-x_j\|}\right)\right)
\\
& = \Hcar^0(x,y) + \mu\ \Hcar^1(x,y),
\end{split}
\end{equation}
where the $0$-th body (the Sun) has mass $1$ and is fixed at the origin and the $i$-th body has mass 
$\mu m_i$ and position-momentum coordinates 
${(x_i,y_i)= (x_{i,1}, x_{i,2}, y_{i,1}, y_{i,2})}$. 
Also, the length and time units are 
chosen so that the gravitational constant is $1$ and the period of an elliptical orbit 
of semi-major axis $1$ is $2\pi$ (so its frequency is $1$, and this the case of the Earth in the Solar system).

The $\mu= 0$ case corresponds to the integrable Keplerian motion of the planets around the Sun 
(no interaction {\em between} planets). Well-known angle-action coordinates for 
the Keplerian motion are Delaunay coordinates.
These are defined body-wise: The 
Delaunay coordinates of the $i$-th body are $( \ell_i, g_i, L_i, G_i) \in \TT^2\times\RR^2$, with $G_i<L_i$ and 
$\TT= \RR/{2\pi\ZZ}$, 
are mapped to the Cartesian coordinates $(x_{i,1}, x_{i,2}, y_{i,1}, y_{i,2})\in \RR^4$ 
 through the following steps:
\[
e_i= \sqrt{1-\left(\frac{G_i}{L_i}\right)^2}, \ a_i=\dfrac{(L_i)^2}{m_i^2}, \ b_i=\dfrac{m_i^2}{L_i},\ E_i = K(\ell_i, e_i),
\]
\[
\begin{pmatrix}
q_{i,1} \\
q_{i,2}
\end{pmatrix}
= a_i 
\begin{pmatrix}
\cos(E_i) - e_i \\
\dfrac{G_i}{L_i}\sin(E_i)
\end{pmatrix},
\qquad
\begin{pmatrix}
x_{i,1} \\
x_{i,2}
\end{pmatrix}
= 
\begin{pmatrix}
 \cos(g_i) & -\sin(g_i)\\
  \sin(g_i) & \cos(g_i)
\end{pmatrix}
\begin{pmatrix}
q_{i,1}\\
q_{i,2}
\end{pmatrix}
\]
\[
\begin{pmatrix}
p_{i,1} \\
p_{i,2}
\end{pmatrix}
=
\dfrac{b_i}{1-e_i \cos(E_i)} 
\begin{pmatrix}
-\sin(E_i) \\
\dfrac{G_i}{L_i}\cos(E_i)
\end{pmatrix}, 
\qquad
\begin{pmatrix}
y_{i,1} \\
y_{i,2}
\end{pmatrix}
= 
\begin{pmatrix}
 \cos(g_i) & -\sin(g_i)\\
  \sin(g_i) & \cos(g_i)
\end{pmatrix}
\begin{pmatrix}
p_{i,1} \\
p_{i,2}
\end{pmatrix}
\]
where $E= K(\ell, e)$ denotes the solution of the Kepler equation $\ell=E-e\sin(E)$.

The Hamiltonian \eqref{eq: hamiltonian cartesian} is then written in Delaunay coordinates $(\ell, g, L, G)$ as 
a function $\Hdel: \TT^{2n}\times\RR^{2n}\to \RR$ given by
\begin{equation*}
\begin{split}
\Hdel(\ell,g,L,G) 
& = \sum_{i=  1}^\Np \frac{- m_i^3}{2 L_i^2} + \mu\ \Hcar^1\comp\Del(\ell,g,L,G) = \Hdel^0(L) + \mu\ \Hdel^1(\ell,g,L,G),\\
\end{split}
\end{equation*}
where $\Del$ denotes the {\em Delaunay map} from Delaunay coordinates $(\ell,g,L,G)$ to Cartesian 
coordinates $(x,y)$ described above.

Let us denote by $\hat G_i = \sum_{1\leq k \leq i} G_k$ the angular momentum of the $i$ first planets. 
It is well-known that the total angular momentum, $\hat G_\Np$, is a first integral of the Hamiltonian system, 
so that we can reduce by one the number of degrees of freedom by fixing the value $\hat G_\Np= \hat G_{\Np,0}$.
By extending the angular momentum map above to a canonical transformation, taking 
$\hat g_i= g_i-g_{i+1}$ for $i= 1,\dots,\Np-1$, and $\hat g_\Np= g_\Np$, 
one gets that $\hat g_\Np$ is a cyclic coordinate in the transformed Hamiltonian in the new coordinates.
Hence, by fixing the total angular momentum $\hat G_\Np$ to a given value $\hat G_{\Np,0}$, 
 one gets a reduced Hamiltonian ${H_{\hat G_{\Np,0}}:\TT^{2n-1}\times\RR^{2\Np-1}\to \RR}$
given by
\begin{equation}\label{eq: hamiltonian reduced}
\Hred_{\hat G_{\Np,0}}(\ell,\hat g, L, \hat G) =  \Hdel^0(L) +  \mu \Hred^1_{\hat G_{\Np,0}}(\ell,\hat g, \L, \hat G),
\end{equation}
with $\hat g= (\hat g_1, \dots, \hat g_{\Np-1})$ and $\hat G= (\hat G_1, \dots, \hat G_{\Np-1})$. 
From now on, we will omit the dependence on $\hat G_{\Np,0}$ from the notation.

A Lagrangian invariant torus of $\Hred$, with a $(2\Np-1)$-dimensional vector of frequencies 
$(\omega^\ell,\omega^{\hat g})\in \RR^\Np\times\RR^{\Np-1}$ and total angular momentum $\hat G_{\Np,0}$, gives raise 
to a Lagrangian invariant torus of $\Hdel$, with a 
$2\Np$-dimensional vector of frequecies $(\omega^\ell, \omega^g)\in \RR^\Np\times\RR^{\Np}$
(see \emph{Reduction Lemma} in \cite{FiguerasHTheorem}). 
The frequencies are related by 
$\omega_i^{\hat g}= \omega_i^g-\omega_{i+1}^g$ for $i= 1, \dots \Np-1$, and $\omega_\Np^{\hat g}= \omega_\Np^g$ is the average of 
$\frac{\partial H}{\partial \hat G_{\Np}}$ over the $(2\Np-1)$-dimensional invariant torus. We emphasize that 
$\omega^\ell$ contains the fast frequencies (the ones coming from the Keplerian motion), 
and that $\omega^{\hat g}$ (or $\omega^g$) contains the slow frequencies
(that in our case are proportional to $\mu$). 
This smallness is a main difficulty when facing the $(1+\Np)$-body problem 
with realistic data (big masses and no big axes).

\subsection{Invariance equation}
In the light of the parameterization method, finding invariant tori for $H$
with frequency vector $\omega= (\omega^\ell, \omega^{\hat g})$ reduces to finding 
a parameterization of the torus 
${K: \TT^{2\Np-1} \to \TT^{2\Np-1} \times \RR^{2\Np-1}}$ satisfying the {\em invariant torus equation}
\begin{equation}
\label{eq:invariance equation}
\Loper K(\theta) +X_H (K(\theta))= 0,
\end{equation}
where $\Lie{\omega}$ is a Lie operator acting on any smooth function $f: \TT^{2\Np-1} \to \RR^M$ by
${\Loper f(\theta) = -{\rm D}f(\theta) \omega,}$
and $X_H= \Omega^{-1}  ({\rm D}H)^\top$ is the Hamiltonian 
vector field with respect to the standard symplectic form given by the matrix
\[
\Omega= \begin{pmatrix} O & -I \\ I & O \end{pmatrix}.
\]
We write $J= \Omega$ when we think of such a matrix as a linear map instead of as a 2-form. 
In particular, we use $J$ to define 
a normal bundle to the torus parameterized by $K$, framed by the columns of 
$N(\theta)= J \Dif K(\theta) (\Dif K(\theta)^\top \Dif K(\theta))^{-1}$, where the columns of 
$\Dif K(\theta)$ frame the tangent bundle.
Moreover, the symmetric matrix 
\[
T(\theta)= N(\theta)^\top \Omega (\Dif X_H(K(\theta)) + J \Dif X_H(K(\theta)) J) N(\theta)
\]
measures how much the normal bundle is twisted (the tangent bundle of an invariant torus 
is fixed). The non-degeneracy of the average of $T$, 
the {\em torsion}, plays the role of the classical Kolmogorov non-degeneracy 
condition in KAM theory.

As it is also costumary in KAM theory, we assume that $\omega$ is Diophantine, i.e. there exists 
$\gamma>0$ and $\tau\geq 2\Np-2$ such that for any $k\in\ZZ^{2\Np-1}$ and $k\neq 0$, 
$|k\cdot\omega|\geq \gamma |k|_1^{-\tau}$. Given any real-analytic function  
$s: \TT^{2\Np-1} \to \RR^M$, we denote by $\Roper(s)$ the only real-analytic function 
$f: \TT^{2\Np-1} \to \RR^M$, with average zero, that satisfies $\Loper f= s-\aver{s}$, 
where $\aver{s}$ denotes the average 
of the function $s$. This operator is in the core of KAM theory.

\section{Continuation from the integrable case with translated tori methods}
\label{section: translated}

Notice that for $\mu=0$ the reduced Hamiltonian \eqref{eq: hamiltonian reduced} 
 has the invariant tori
\begin{equation}
\label{def:plane torus}
K_{\hat G_0}(\theta^\ell, \theta^{\hat g})= (\theta^\ell, \theta^{\hat g}, L_0, {\hat G_0})
\end{equation}
where the components of $L_0$ are determined by the masses of the bodies and the fast frequencies $\omega^\ell$ 
(by the third Kepler's law), but the secular frequency $\omega^{\hat g}_0$ is zero 
(not $\omega^{\hat g}$!) and $\hat G_0$ is free: 
there is an $(\Np-1)-$parameter family of $(2\Np-1)$-dimensional tori foliated by 
$\Np$-dimensional invariant tori. As a result, 
the torsion is noninvertible (since there is no twist in the $\hat G$ direction). 
In summary: the problem is degenerate.

\subsection{A translated torus algorithm}

As mentioned above, the degeneracy of the problem imposes a first obstacle
for applying any numerical KAM scheme for performing any continuation with respect to $\mu$.
In the spirit of \cite{GonzalezHL13, GonzalezHL22}, we can overcome this degeneracy 
by introducing a counterterm $\lambda \, \Pi_{\hat G}$ to  the Hamiltonian
\eqref{eq: hamiltonian reduced}, where $\Pi_{\hat G}:\TT^{2n-1}\times\RR^{2n-1}\to \RR$ 
is the projection onto
the $\hat G$ coordinate. Hence, by denoting $X_{\hat G}= X_{\Pi_{\hat G}}$,
instead of solving Equation 
\eqref{eq:invariance equation}, we solve the extended system 
\begin{equation}
\label{eq:translated}
\left\{
\begin{split}
\Lie{\omega} K(\theta)+X_H(K(\theta))+ X_{\hat G}(K(\theta)) \lambda & = 0,\\
\aver{\Pi_{\hat G}(K(\theta))} -  \hat G_0 & = 0 ,
\end{split}
\right.
\end{equation}
for a fixed constant $\hat G_0$.
The {\em invariant tori} satisfying \eqref{eq:translated} are {\em translated tori} 
for the original Hamiltonian system. 
Since $\hat G_0$ is an extra parameter, under appropriate non-degeneracy conditions 
(that we will see later are very mild), 
we can find families of translated tori, labeled by $\hat G_0$. The use of 
counterterms in KAM theory goes back to 
the works of Moser and Herman \cite{Moser67,Herman98,Fejoz04,Fejoz16}.

Notice that,  for a given $\hat G_0$, for $\mu= 0$ the parameterization \eqref{def:plane torus} 
satisfies \eqref{eq:translated} with frequency 
$\omega= (\omega^\ell,\omega^{\hat g})$, by selecting $\lambda= \omega^{\hat g}$. 
The idea is then performing a continuation method for solving the {\em translated torus equation} 
\eqref{eq:translated}  for couples $(K,\lambda)$ up to the value $\mu= \mu_0$. The rationale behing 
this method is that, if there were an invariant torus with frequency 
$\omega= (\omega^\ell,\omega^{\hat g})$ and 
$\aver{\Pi_{\hat G}\comp K} =  \hat G_0$ for $\mu= \mu_0$, then after the continuation procedure we would find 
$(K,\lambda)$ with  $\lambda= 0$. Using perturbation theory up to order one 
(expanding in Poincaré-Lindsedt series)
we get an approximation of 
$\hat G_0$ by solving the equation
$\aver{ \Pi_{\hat g} X_{H^1}\circ K_{\hat G}} = \omega^{\hat g}/\mu_0,$
where $\Pi_{\hat g}$ is the projection onto the $\hat g$ component. We then 
continue this solution from $\mu=0$ to $\mu=\mu_0$ by solving the equations at each step, 
see Subsection \ref{section: solution translation}.
Finally, 
since this value of $\hat G_0$ is not exact, we don't get $\lambda=0$ at $\mu=\mu_0$. However, 
at $\mu_0$ on can also tune $\hat G_0$ to get $\lambda = 0$ by a Newton method (again using 
perturbation methods for computing $\frac{\partial \lambda}{\partial \hat G_0}$).
\footnote{In applications, estimates of $\hat G_0$ could also be obtained by 
methods such as averaging the $\hat G$
components of a quasiperiodic orbit obtained using frequency analysis \cite{Laskar05,GomezMondeloSimo10,LuqueV09,DasSaikiSanderYorke17}}.

\subsubsection{Solving Equations \eqref{eq:translated}}
\label{section: solution translation}
More concretetly, from an approximate solution $(K(\theta), \lambda)$ of \eqref{eq:translated} we can perform 
a quasi-Newton correction 
of the form $(P(\theta)\xi(\theta), \Delta\lambda)$, where the matrix 
$P(\theta)= \begin{pmatrix} \Dif K(\theta) & N(\theta))\end{pmatrix}$,
obtained by yuxtaposing the tangent and normal frames given above, is approximately symplectic.
Taking into account that 
the inverse of $P(\theta)$ is close to $-\Omega P(\theta)^T\Omega$, 
we end up with the linear system 
\[
\left\{
\begin{split}
\Loper \xiL(\theta)+ T(\theta)\xiN(\theta)+ \bL(\theta) \Delta\lambda  &= \etaL(\theta),\\
\Loper \xiN(\theta)+\phantom{T(\theta)\xiN(\theta) +}  b^N(\theta) \Delta\lambda &= \eta^N(\theta),\\
\aver{ \Pi_{\hat G} (\Dif K(\theta)\xiL(\theta)+N(\theta)\xiN(\theta))} & = \eta^{\hat G_0},
\end{split}
\right.
\]
where  
$\begin{pmatrix}\bL(\theta)\\ \bN(\theta)\end{pmatrix}= 
\begin{pmatrix}N(\theta)^\top\\-\DK(\theta)^\top\end{pmatrix}\Omega X_{\hat G}(K(\theta))$, 
$\begin{pmatrix}\etaL(\theta)\\\etaN(\theta)\end{pmatrix}=
\begin{pmatrix}-N(\theta)^\top\\ \DK(\theta)^\top\end{pmatrix}\Omega (\Lie{\omega} K(\theta)+X_H(K(\theta))+\lambda X_{\hat G}(K(\theta)))$, 
$\eta^{\hat G_0}= -\aver{\Pi_{\hat G}\circ K} + \hat G_0$.

As it is customary in these types of schemes, one solves them up to some a-priori unknowns:
In this case
the average $\xiN_0= \aver{\xiN}$ and $\Delta\lambda$. These two satisfy the linear system 
\begin{equation}\label{eq: averages}
\begin{pmatrix}
\aver{T}  & \aver{\tbL} \\
\aver{\Pi_{\hat G} (N -DK \Roper T)} 
& - \aver{\Pi_{\hat G} ( DK \Roper\tbL + N \Roper b^N)}
\end{pmatrix}
\begin{pmatrix}
\xiN_0 \\
\Delta\lambda
\end{pmatrix}
=
\begin{pmatrix}
\aver{\tetaL} \\ \eta^{\hat G_0} - \aver{\Pi_{\hat G} \Dif K \Roper\tetaL }
\end{pmatrix},
\end{equation}
where $\tbL= \bL-T \Roper b^N$ and $\tetaL= \etaL-T\Roper \eta^N$.
If the matrix in \eqref{eq: averages}, to which we will refer to as the {\em supertorsion} $\aver{\hat T}$, 
is regular, then the method can continue by computing
\[
\left\{
\begin{split}
\xiN(\theta) & = \xiN_0+ \Roper\eta^N(\theta) - \Roper b^N(\theta) \Delta\lambda,\\
\xi^L(\theta) & =  \Roper\tetaL(\theta) - 
\Roper T(\theta) \xiN_0 - \Roper\tbL(\theta) \Delta\lambda,
\end{split}
\right.
\]
and at the next step 
we get a quadratically better estimate $(K+P\xi , \lambda+ \Delta\lambda)$.

\begin{remark}
In particular, in the case $\mu= 0$,  the torsion and the supertorsion of the torus \eqref{def:plane torus}
are
\[
\aver{T} = 
\begin{pmatrix}
\Dif^2 H^0(L_0) & O \\
O & O  
\end{pmatrix}, \quad 
\aver{\hat T} = 
\begin{pmatrix}
\Dif^2 H^0(L_0) & O & O \\
O & O  & I                       \\
O & I & O
\end{pmatrix},
\]
respectively.
Notice that, even though the torsion 
is degenerate, the supertorsion is not, permiting to start the continuation of {\em translated tori} from $\mu= 0$.

\end{remark}

\subsection{An invariant torus algorithm}

Once the continuation explained in Subsection \ref{section: translated}
reaches the parameter value $\mu_0$, one obtains a translated torus $K$ with 
translation $\lambda$ that, ideally, should be zero. In order to refine the (approximate) invariant torus, we follow 
a similar scheme as deviced before but for Equation \eqref{eq:invariance equation}, 
in which the only unknown is $K$. (Another possibility is following the previous scheme, and 
tune parameter $\hat G$ so that one gets $\lambda= 0$.)
Then, given an approximate solution $K$ of \eqref{eq:invariance equation}, its correction is given by 
$P(\theta) \xi(\theta)$. After truncating 
up to linear terms we obtain the linear system 
\[
\left\{
\begin{split}
\Loper \xiL(\theta)+ T(\theta)\xiN(\theta) &= \etaL(\theta), \\
\Loper \xiN(\theta)+\phantom{T(\theta)\xiN(\theta)}  &= \eta^N(\theta),
\end{split}
\right.
\]
where $\begin{pmatrix}\etaL(\theta)\\\etaN(\theta)\end{pmatrix}
= \begin{pmatrix}-N(\theta)^\top\\ \DK(\theta)^\top\end{pmatrix} \Omega (\Lie{\omega} K(\theta)+X_H(K(\theta))$. 

Notice that in this case we assume that the torsion $\aver{T}$ is non-degenerate, so that this last system 
can be solved as
\[
\left\{
\begin{split}
\xiN_0  &=  \aver{T(\theta)}^{-1} \aver{\etaL(\theta) -  T(\theta)\Roper\eta^N(\theta)}, \\
\xiN(\theta) & = \xiN_0+ \Roper\eta^N(\theta),\\
\xiL(\theta) & =  \Roper\left(\etaL(\theta)-T(\theta)\xiN(\theta)\right).
\end{split}
\right.
\]

\section{Application to the Sun-Jupiter-Saturn problem}

We have implemented the algorithms discussed in Section \ref{section: translated} 
in C++ 
(see e.g. \cite{HaroCFLM16} for similar implementations). 
A key point of the implementation is to use FFT routines for fast evaluations of the vector 
field on the parameterization, and for evaluating the operator $\Roper$. 
We have adapted the FFT routines from \cite{PressTVF02}
to work with multiprecision arithmetics with \text{mpfr} (see \cite{mpfrlib}).
Another key point is parallelization using {\tt openmp} (see \cite{book_openmp}).  

Here we present the specifics for the planar Sun-Jupiter-Saturn problem with realistic 
values of their parameters(masses, frequencies, ephemerides...)
The source for the values of the parameters we have used come from astronomical observations
from NASA: \cite{planetarysheet, sunsheet}.
(Other important tools such as frequency analysis 
\cite{Laskar05,GomezMondeloSimo10,LuqueV09,DasSaikiSanderYorke17} could have also been used 
for getting these data.)

The masses of Jupiter and Saturn are $0.9546\cdot 10^{-3}$ and $0.2856\cdot 10^{-3}$, respectively, 
so $m_1= 0.9546$, $m_2= 0.2856$ and $\mu_0= 10^{-3}$.
From their orbital elements and Kepler laws (hence using semiaxis $a_1,a_2$ and 
excentricities $e_1,e_2$ of Jupiter and Saturn, respectively) 
we get the approximations for the frequencies of the Keplerian motions, 
say
\begin{equation}
\label{eq:omegaell}
\omega^\ell=  (8.39549288702546301204\cdot 10^{-2}, 3.38240117059304358259\cdot 10^{-2}).
\end{equation}
From the ephemerides, and in particular their precession motions, we get 
\begin{equation}
\label{eq:omegahat}
\omega^{\hat g_1}= \omega^{g_1}-\omega^{g_2} = -1.85007988077595000000\cdot 10^{-5}.
\end{equation}
Also, the total angular momentum is approximately
\[
\hat G_{2,0}= 3.05839852910896096675.
\]

While the total angular momentum $\hat G_2$ is preserved (and equal to the value 
$\hat G_{2,0}$), the one of Jupiter, 
$\hat G$, is not. 
An approximation of the average of the angular momentum of Jupiter 
comes from the Kepler approximation, and approximations of order $\mu_0$ are obtained by solving the
equation $\aver{ \Pi_{\hat g} X_{H^1}\circ K_{\hat G}} = \omega^{\hat g}/\mu_0.$
After several iterations of the secant method we obtain the approximation
$
\hat G_0= 2.17647359010273488684.
$
These  preliminary computations are performed with {\tt long double} precision C 
arithmetic and the parameterizations are given by 
grids of size $128^3$.

\begin{figure}
\begin{tabular}{ccc}
\includegraphics[width=0.3\linewidth]{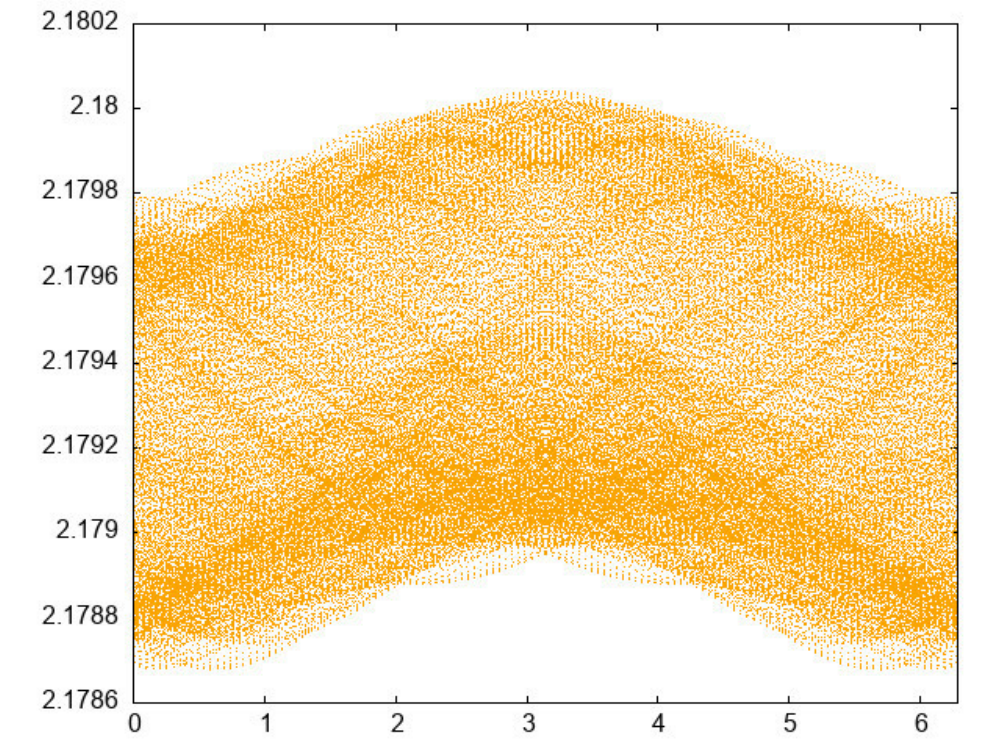} &
\hspace{-0.5cm}
\includegraphics[width=0.3\linewidth]{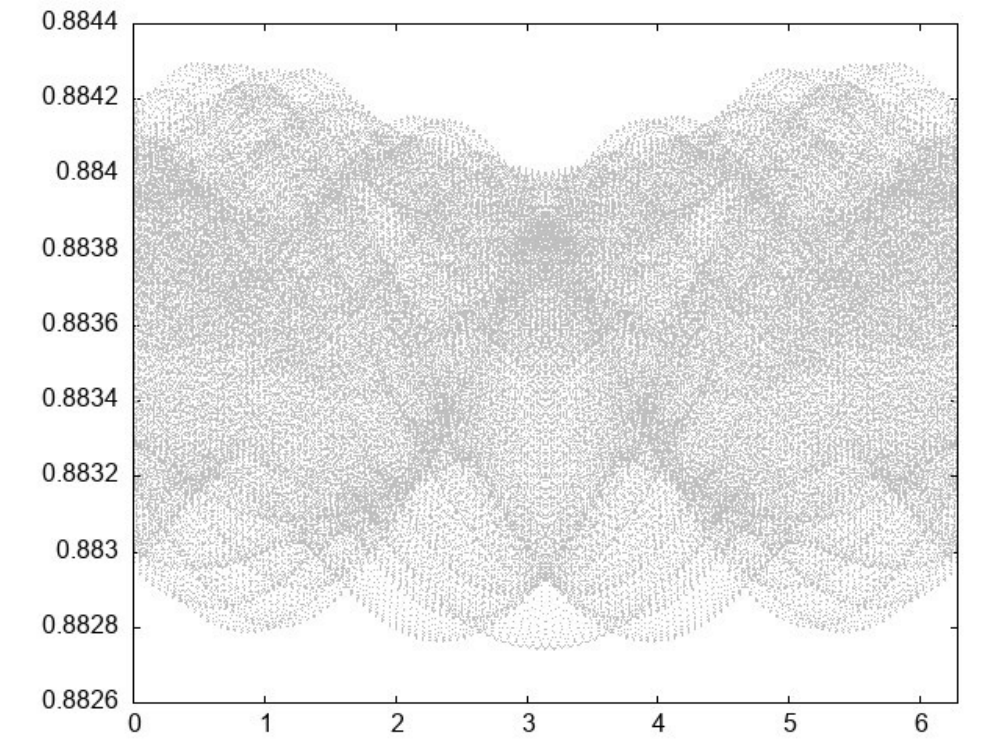} &
\hspace{-0.5cm}
 \includegraphics[width=0.3\linewidth]{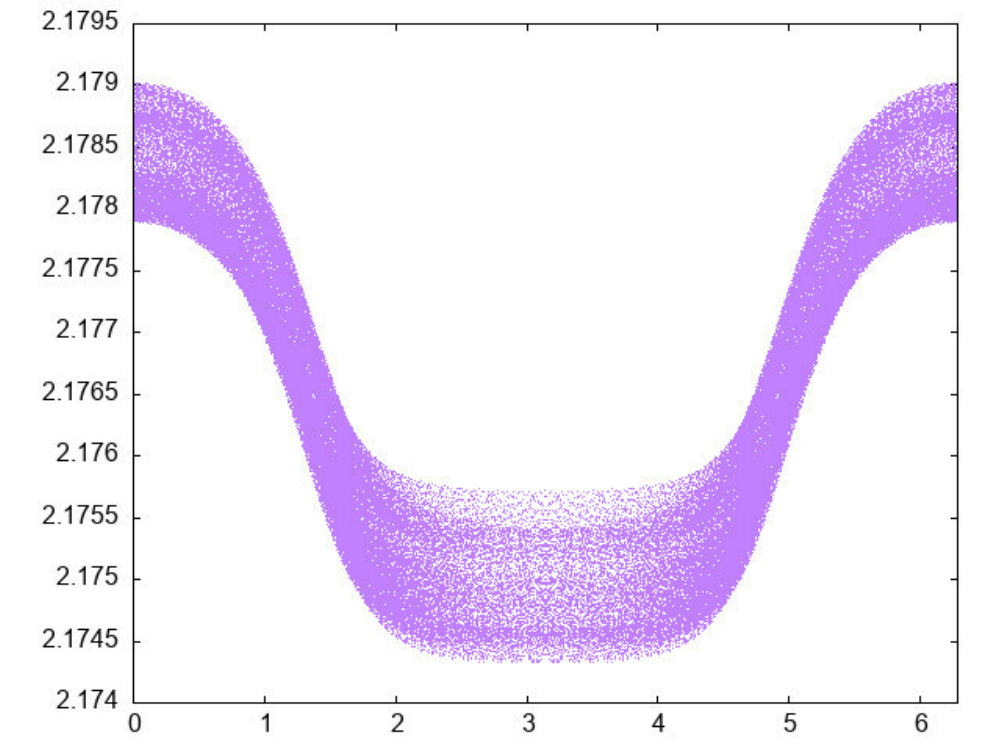} 
\end{tabular}
\caption{\label{fig:delaunay-projection} Projections of the 3D invariant torus in  Delaunay coordinates 
onto $(\ell_1,L_1)$, $(\ell_2,L_2)$ and $(\hat g,\hat G)$ components.}
\end{figure}

\begin{figure}
\begin{tabular}{cc}
\includegraphics[width=0.4\linewidth]{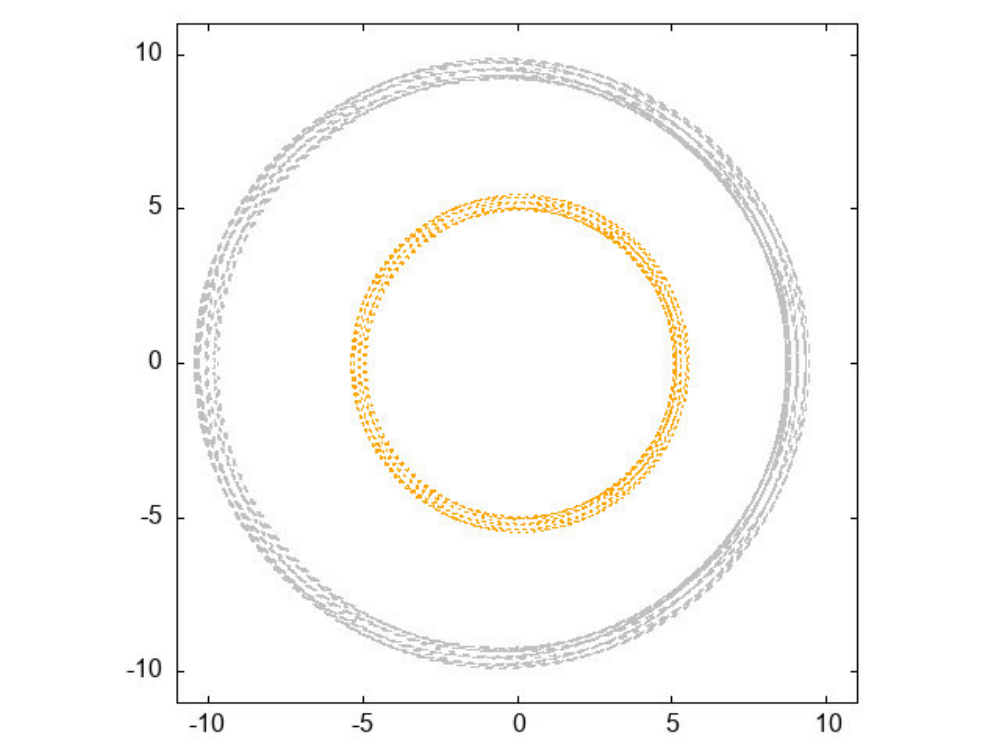} &
\includegraphics[width=0.4\linewidth]{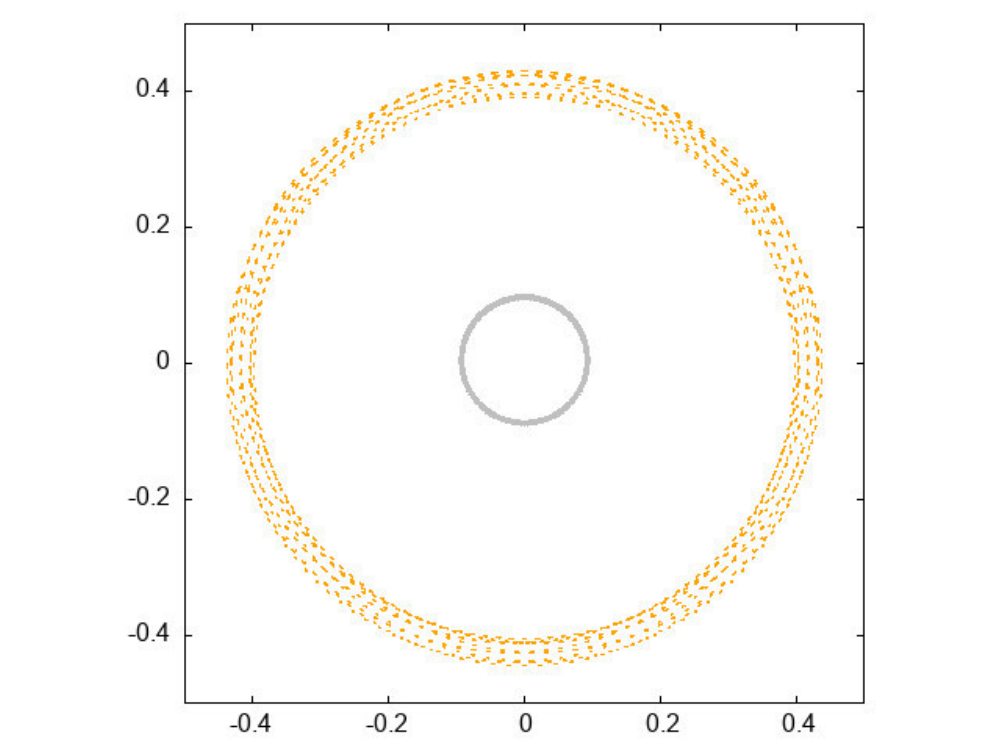} 
\end{tabular}
\caption{\label{fig:cartesian-projection} Projections of the 3D invariant torus (generating a 4D torus) in Cartesian coordinates 
onto positions $x_1= (x_{1,1},x_{1,2})$, 
$x_2= (x_{2,1},x_{2,2})$ and momenta $y_1= (y_{1,1},y_{1,2})$, 
$y_2= (y_{2,1},y_{2,2})$. 
Coordinates $x_1,y_1$ correspond to Jupiter and $x_2,y_2$ correspond to Saturn,  and are plot in orange and grey, respectively.
}
\end{figure}

The first thing we have done is finding the value of $\hat G_0$ such that 
at $\mu=\mu_0$, when performing the continuation of the translated torus algorithm, 
we obtain $\lambda=0$.
To do this we performed the following 4 times: 
At $\mu=0$ we have the invariant torus \eqref{def:plane torus} with an approximation of 
the desired $\hat G_0$ value. Then 
we perform a continuation using the translated 
torus algorithm from $\mu= 0$ to $\mu= 10^{-3}$ and get 
an approximately translated torus (that, unfortunately, does not satisfy $ \lambda=0$).
Finally, by doing a Newton step for solving  $\lambda(\hat G_0)= 0$ we obtain a better estimate 
of $\hat G_0$. Then we repeat the process. By doing this, in the first run we obtain 
an approximately translated torus with invariance error $4.8\cdot 10^{-9}$, 
moment error $ \aver{\Pi_{\hat G}K}-\hat G_0=1.5\cdot 10^{-14}$, and 
$\lambda= -8.67345532598763273085\cdot 10^{-7}$. The Newton step gives us 
a better estimate $\hat G_0= 2.17658253666877214401$. After the fourth time we do this, we obtain 
an approximately translated torus with 
$\hat G_0= 2.17657425006565519231$, invariance error 
$5.3\cdot 10^{-10}$, moment error $\aver{\Pi_{\hat G}K}-\hat G_0= 6.5\cdot 10^{-17}$ and, 
$\lambda= 3.14309229785154830993\cdot 10^{-14}$. For this last torus the 
supertorsion \eqref{eq: averages} is 
\[
\left(
\begin{tabular}{rrr|r}
$-1.15830235 \cdot 10^{-1}$ &	$1.43198239\cdot 10^{-3}$ &	
$7.41752484\cdot 10^{-4}$ &	$2.11280898\cdot 10^{-1}$ \\
$1.43198239\cdot 10^{-3}$  & $-1.23922829\cdot 10^{-1}$ &	
$-5.21182674\cdot 10^{-3}$ &	$-3.45764287\cdot 10^{-1}$ \\
$7.41752484\cdot 10^{-4}$  & $-5.21182674\cdot 10^{-3}$ &	
$-3.18265383\cdot 10^{-3}$ &	$5.20444218\cdot 10^{-1}$ \\
       \hline
$2.11280898\cdot 10^{-1}$ & $-3.45764287\cdot10^{-1}$ &	
$5.20444218\cdot 10^{-1}$ &	$8.38171177\cdot 10^0\,\,\,\,$
\end{tabular}
\right)
\]
The norm of the inverse of the supertorsion is  $3.9\cdot 10^{1}$
and of the inverse of the torsion is
$3.5\cdot 10^{2}$.
The whole computation takes less than one hour, 
with the first continuation taking around 17 minutes, and the last around 12 minutes.
Different projections of the invariant torus are shown in 
Figures~\ref{fig:delaunay-projection} and \ref{fig:cartesian-projection}.

Later on we refined the approximately invariant torus using the invariant torus algorithm
by increasing the 
precision with 
{\tt long double}, {\tt \_\_float128} and, finally {\tt mpfr}. We gradually increased the accuracy and the 
size of the grids. In the last run, 
the input torus was given with a grid of size $512^3$ with 57 digits, and the output with 
a grid of size $1024^3$ 
with 76 digits. The input error was $9.1\cdot 10^{-29}$, and the error saturated at the first step 
to $3.9\cdot 10^{-54}$ (the error at the second step was $2.9\cdot 10^{-54}$). Although it is 
not used during the computations, we obtain 
that this torus has $\hat G_0= 2.17657418883872689352084685277943$.
In this case, one Newton step took around one week, and the top size of RAM memory used was 194G. 

\begin{figure}
\begin{tabular}{cc}
$\theta^{\ell_1}$ \includegraphics[width=0.4\linewidth]{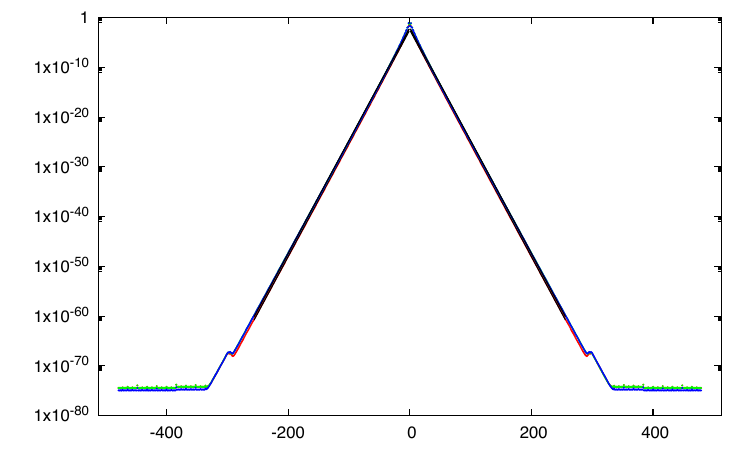} &
$\theta^{\ell_2}$ \includegraphics[width=0.4\linewidth]{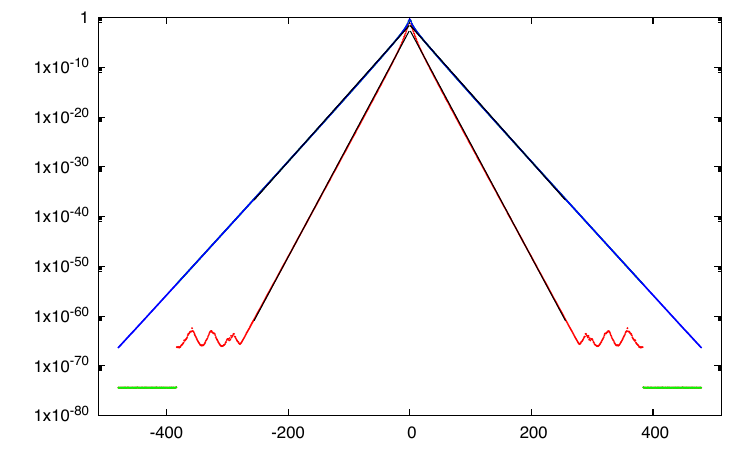} \\
$\theta^{\hat g_1}$ \includegraphics[width=0.4\linewidth]{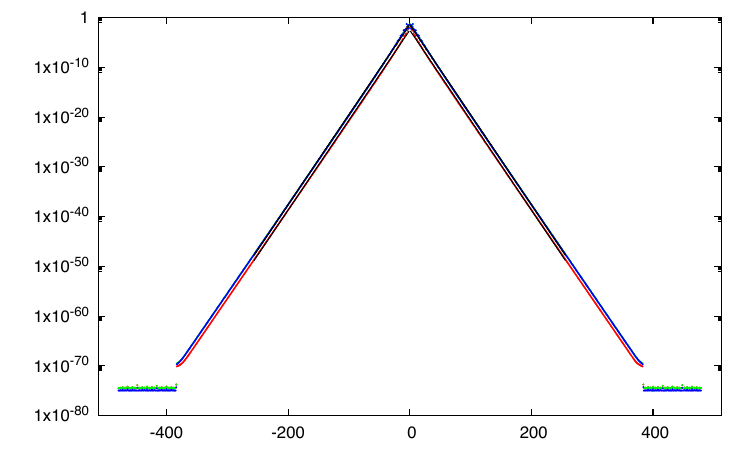} &
{
\begin{minipage}{0.4\linewidth}
\vspace{-3cm}
\begin{tabular}{|c|ccc|}
\hline 
& {\footnotesize $K^{\ell_1} + {\ii} K^{L_1}$}
& {\footnotesize $K^{\ell_2} + {\ii} K^{L_2}$}  
& {\footnotesize $K^{\hat g} + {\ii} K^{\hat G}$}   \\
\hline
$\rho_{\ell_1}$          & 0.524453%710213535%  
& 0.521794%003886631% 
& 0.522239%216823448% 
\\
$\rho_{\ell_2}$          & 0.416514%34944186%    
&0.417078%447086004% 
& 0.417001%136956824%
\\
$\rho_{\hat g}$      & 0.525773%529859222%  
& 0.316526%890290233% 
& 0.316557%291213263%
\\
\hline
\end{tabular}
\end{minipage}
}
\end{tabular}
\caption{
\label{fig:fourier}
Fits of Fourier coefficients of the (complexified) components of the parameterization
$\tilde f_{\ell_1}, \tilde f_{\ell_2}$ and $\tilde f_{\hat g}$, 
and estimates of analyticity strips.}
\end{figure}

From this last torus we have estimated how fast the Fourier coefficients decrease and, so, 
its analiticity radius. To do so, we have 
fit the Fourier coefficients of the complexifications $K^{\ell_1} + {\ii} K^{L_1}$, 
$K^{\ell_2} + {\ii} K^{L_2}$ and  $K^{\hat g} + {\ii} K^{\hat G}$ with respect to
each of the angles $\theta^{\ell_1}, \theta^{\ell_2}, \theta^{\hat g}$, thus obtaining estimates of the 
 the analyticity strips $\rho_{\ell_1}, \rho_{\ell_2}, \rho_{\hat g}$.
The results are shown in 
Figure~\ref{fig:fourier} where, for a Fourier expansion 
\[
f(\theta^{\ell_1},\theta^{\ell_2}, \theta^{\hat g}) = 
\sum_{k_{\ell_1},k_{\ell_2}, k_{\hat g}} f_{k_{\ell_1},k_{\ell_2}, k_{\hat g}} 
e^{\ii (k_{\ell_1} \theta^{\ell_1}+ k_{\ell_2} \theta^{\ell_2}+ k_{\hat g} \theta^{\hat g})},
\]
we fit the analyticity strips $\rho_{\ell_1}, \rho_{\ell_2}, \rho_{\hat g}$ of each of the angles
by considering the univariate Fourier series  
\[
\tilde f_{\ell_1}(\theta^{\ell_1})= \sum_{k_{\ell_1}} \left(\sum_{k_{\ell_2}, k_{\hat g}} |f_{k_{\ell_1},k_{\ell_2}, k_{\hat g}}|\right)
e^{\ii (k_{\ell_1} \theta^{\ell_1})}, \quad
\tilde f_{\ell_2}(\theta^{\ell_2})= \sum_{k_{\ell_2}} \left(\sum_{k_{\ell_1}, k_{\hat g}} |f_{k_{\ell_1},k_{\ell_2}, k_{\hat g}}|\right)
e^{\ii (k_{\ell_2} \theta^{\ell_2})},
\]
\[
\tilde f_{\hat g}(\theta^{\hat g})= \sum_{k_{\hat g}} \left(\sum_{k_{\ell_1}, k_{\ell_2}} |f_{k_{\ell_1},k_{\ell_2}, k_{\hat g}}|\right)
e^{\ii (k_{\hat g} \theta^{\hat g})}
\]
and doing a standard fit on their coefficients.

\section{Numerical verification of the KAM constants}

The numerical certification of the existence of the invariant torus is based on the KAM Theorem 
and the Iterative Lemma appearing in \cite{FiguerasHTheorem}. For the sake of completeness, we include their tailored and simplified versions (with the most relevant
hypotheses) in appendix
\ref{section: theorems},
 so it will guide us in all the data needed for doing the validation. For the specific 
expression of all the constants we refer the reader to \cite{FiguerasHTheorem}, where 
they appear in the appendices.

Given the $\omega=(\omega^\ell, \omega^{\hat g_1})$ in \eqref{eq:omegaell}, \eqref{eq:omegahat} 
we can certify (using the validation techniques in \cite{FiguerasHL17}) that at distance $10^{-80}$ there is 
a Diophantine vector with $\tau = 2.4$ and $\gamma = 1.69 \cdot 10^{-6}$. Moreover, we choose 
the radius of analiticity to be
$\rho = 0.1$ and $\delta=\frac{\rho}{6}$.

The hypotheses in $H_1$ control the Hamiltonian and its associated vector field in a tubular 
neighborhood of the torus $K(\TT^m_{\rho})$. In our case, it is enough to take these constants 
to be
\[
\cteXH = 0.09,
\cteDXH = 129,
\cteDXHT = 129,
\cteDDXH = 5\cdot 10^{10}.
\]

The hypotheses in $H_2$ control the parameterization $K$ and all the geometric infomation 
it has (the bunbles $DK$, $N$ and so on). In our case, these constants are computed 
with the approximation and obtained
\[
\|DK\| = 4.7811815833,
\|DK^T\| = 6.8755882886,
\|B\| = 7.35806411265,
\]
\[
\|\N\|=  3.5704498717, 
\|\NT\|=  2.8621242724, 
|\aver{T}|^{-1}= 354.07743243.
\]
The corresponding $\sigma$ constants are obtained by multipling these norms by a factor 
$1+10^{-10}$.

With 
this information
we can run the Iterative Lemma several steps, say 10,
and with different initial invariance errors 
and 
then apply the KAM Theorem to see if it converges (Inequality \eqref{eq:KAM:HYP} is fulfilled). We have obtained 
that with $\|\etaL\|=10^{-38}, \|\etaN\|=10^{-44}$. However, from our numerics we obtain 
that our torus satisfies
$\|\etaL\|= 3.6\cdot 10^{-54}$ and 
$\|\etaN\|= 1.7\cdot 10^{-57}$, which are very much smaller than the thresholds!

\section{Computation details}
For running the continuation method on the translated torus algorithm
from the integrable system, we 
run the programs is an out to date MacBook Air laptop
with one CPU 1.7 GHz Dual-Core Intel i7 and RAM memory 8G, since for the approximation we 
work with {\tt long double} C arithmetics and the tori are discretized in $128^3$ 
nodes, accounting to 32M of memory for each of the six  
components of the parameterization of the torus. We have also adapted and tested the 
programs to work with quadruple precision {\tt \_\_float128} C arithmetics. 
For the invariant torus algorithm, we have used an iMac Pro with one CPU 3,2 GHz Intel Xeon W 
with 8 cores and RAM Memory 256G, working
with several extended precision arithmetics with 
{\tt mpfr} (up to 76 decimal digits, that correspond to 64 bytes, respectively) and the torus is 
discretized in $1024^3$ nodes 
, accounting 64G of memory
for each to the components.
This last computation has also been run in the UPPMAX supercomputer.

Finally, we give some numbers to provide an idea of the order of magnitude of the managed data structures
at the final stages of the computations.
The data structures are 
complex vectors, that store couples of real grids. Moreover, handling of memory (both RAM and disc) by 
{\tt mpfr} is anisotropic. For instance, for the computation of the torus with 
76 digits the program uses up to  194G of RAM memory for 
handling one single complex grid of size $1024^3$, and 891G of memory  disk to store the objects being computed
by the program. For files storing the same number of {\tt mpfr} objects, $2\cdot1024^3/8$,  the sizes range from 87M to 8.8G.
 
\section{Acknowledgements}

The authors are grateful to Alejandro Luque, Kristian Bjerklöv, Chiara Caracciolo and Andreas Strömbergsson 
for fruitful discussions. 

J.-Ll.F. has been partially supported by the Swedish VR Grant 2019-04591, and 
A.H. has  been supported by the Spanish grant PID2021-125535NB-I00 (MCIU/AEI/FEDER, UE), 
and by the Spanish State Research Agency, through the
Severo Ochoa and Mar\'ia de Maeztu Program for Centers and Units of
Excellence in R\&D (CEX2020-001084-M).
Some computations were enabled by resources in project 
NAISS 2023/5-192 provided by the National Academic Infrastructure for Supercomputing in Sweden (NAISS) 
at UPPMAX, funded by the Swedish Research Council through grant agreement no. 2022-06725.

\bibliographystyle{plain}
\bibliography{references}

\def\cprime{$'$} \def\cprime{$'$} \def\cprime{$'$} \def\cprime{$'$}
\begin{thebibliography}{10}

\bibitem{planetarysheet}
Planetary fact sheet - metric.
\newblock \url{https://nssdc.gsfc.nasa.gov/planetary/factsheet/}.
\newblock Accessed: 2023-12-05.

\bibitem{sunsheet}
Sun fact sheet.
\newblock \url{https://nssdc.gsfc.nasa.gov/planetary/factsheet/sunfact.html}.
\newblock Accessed: 2023-12-05.

\bibitem{Arnold62}
V.~I. Arnold.
\newblock On the classical perturbation theory and the stability problem of
  planetary systems.
\newblock {\em Dokl. Akad. Nauk SSSR}, 145:487--490, 1962.

\bibitem{Arnold63a}
V.I. Arnold.
\newblock Proof of a theorem of {A}. {N}. {K}olmogorov on the preservation of
  conditionally periodic motions under a small perturbation of the
  {H}amiltonian.
\newblock {\em Uspehi Mat. Nauk}, 18(5 (113)):13--40, 1963.

\bibitem{Arnold63b}
V.I. Arnold.
\newblock Small denominators and problems of stability of motion in classical
  and celestial mechanics.
\newblock {\em Russ. Math. Surveys}, 18:85--192, 1963.

\bibitem{BroerHS96}
H.W. Broer, G.B. Huitema, and M.B. Sevryuk.
\newblock {\em Quasi-periodic motions in families of dynamical systems. {O}rder
  amidst chaos}.
\newblock Lecture Notes in Math., {V}ol 1645. Springer-Verlag, Berlin, 1996.

\bibitem{CallejaCLa}
R.~Calleja, A.~Celletti, and R.~de~la Llave.
\newblock A {KAM} theory for conformally symplectic systems: efficient
  algorithms and their validation.
\newblock {\em J. Differential Equations}, 255(5):978--1049, 2013.

\bibitem{Castan2017thesis}
T.~Castan, J.~F{\'e}joz, A.~Chenciner, L.N. ), A.I. Neishtadt, L.C. ), J.P.M.
  ), V.Y. Kaloshin, E.~S{\'e}r{\'e}, Hauts-de-Seine / 1992-....). {\'E}cole
  doctorale Astronomie et astrophysique d'{\^I}le-de France~(Meudon, et~al.
\newblock {\em Stability in the Plane Planetary Three-body Problem}.
\newblock 2017.

\bibitem{CellettiC95}
A.~Celletti and L.~Chierchia.
\newblock A constructive theory of {L}agrangian tori and computer-assisted
  applications.
\newblock In {\em Dynamics Reported}, pages 60--129. Springer, Berlin, 1995.

\bibitem{book_openmp}
R.~Chandra, L.~Dagum, D.~Kohr, R.~Menon, D.~Maydan, and J.~McDonald.
\newblock {\em Parallel programming in OpenMP}.
\newblock Morgan kaufmann, 2001.

\bibitem{Chierchia03}
L.~Chierchia.
\newblock K{AM} lectures.
\newblock In {\em Dynamical systems. {P}art {I}}, Pubbl. Cent. Ric. Mat. Ennio
  Giorgi, pages 1--55. Scuola Norm. Sup., Pisa, 2003.

\bibitem{ChierchiaP11}
L.~Chierchia and G.~Pinzari.
\newblock The planetary {$N$}-body problem: symplectic foliation, reductions
  and invariant tori.
\newblock {\em Invent. Math.}, 186(1):1--77, 2011.

\bibitem{DasSaikiSanderYorke17}
S.~Das, Y.~Saiki, E.~Sander, and J.~A. Yorke.
\newblock Quantitative quasiperiodicity.
\newblock {\em Nonlinearity}, 30(11):4111--4140, 2017.

\bibitem{Llave01}
R.~de~la Llave.
\newblock A tutorial on {KAM} theory.
\newblock In {\em Smooth ergodic theory and its applications (Seattle, WA,
  1999)}, volume~69 of {\em Proc. Sympos. Pure Math.}, pages 175--292. Amer.
  Math. Soc., Providence, RI, 2001.

\bibitem{GonzalezJLV05}
R.~de~la Llave, A.~Gonz{\'a}lez, {\`A}.~Jorba, and J.~Villanueva.
\newblock K{AM} theory without action-angle variables.
\newblock {\em Nonlinearity}, 18(2):855--895, 2005.

\bibitem{Dumas14}
H.S. Dumas.
\newblock {\em The {KAM} story}.
\newblock World Scientific Publishing Co. Pte. Ltd., Hackensack, NJ, 2014.
\newblock A friendly introduction to the content, history, and significance of
  classical Kolmogorov-Arnold-Moser theory.

\bibitem{Fejoz04}
J.~F{\'e}joz.
\newblock D\'emonstration du `th\'eor\`eme d'{A}rnold' sur la stabilit\'e du
  syst\`eme plan\'etaire (d'apr\`es {H}erman).
\newblock {\em Ergodic Theory Dynam. Systems}, 24(5):1521--1582, 2004.

\bibitem{Fejoz16}
Jacques F\'{e}joz.
\newblock Introduction to {KAM} theory with a view to celestial mechanics.
\newblock In {\em Variational methods}, volume~18 of {\em Radon Ser. Comput.
  Appl. Math.}, pages 387--433. De Gruyter, Berlin, 2017.

\bibitem{FiguerasHValidation}
J.-Ll Figueras and A.~Haro.
\newblock {A Computer-Assisted Proof of the Existence of Quasiperiodic
  Solutions of the planar Sun-Jupiter-Saturn Problem with Realistic Data}.
\newblock {\em (In progress)}.

\bibitem{FiguerasHTheorem}
J.-Ll Figueras and A.~Haro.
\newblock {A modified parameterization method for invariant Lagrangian tori for
  partially integrable Hamiltonian systems}.
\newblock {\em (Accepted in Physica D)}.

\bibitem{FiguerasHL17}
J.-Ll. Figueras, A.~Haro, and A.~Luque.
\newblock Rigorous {C}omputer-{A}ssisted {A}pplication of {KAM} {T}heory: {A}
  {M}odern {A}pproach.
\newblock {\em Found. Comput. Math.}, 17(5):1123--1193, 2017.

\bibitem{FiguerasHL20}
J.-Ll. Figueras, A.~Haro, and A.~Luque.
\newblock Effective bounds for the measure of rotations.
\newblock {\em Nonlinearity}, 33(2):700--741, 2020.

\bibitem{mpfrlib}
L.~Fousse, G.~Hanrot, V.~Lef\`{e}vre, P.~P\'{e}lissier, and P.~Zimmermann.
\newblock Mpfr: A multiple-precision binary floating-point library with correct
  rounding.
\newblock {\em ACM Trans. Math. Softw.}, 33(2):13–es, jun 2007.

\bibitem{GomezMondeloSimo10}
G.~G\'{o}mez, J.-M. Mondelo, and C.~Sim\'{o}.
\newblock A collocation method for the numerical {F}ourier analysis of
  quasi-periodic functions. {I}. {N}umerical tests and examples.
\newblock {\em Discrete Contin. Dyn. Syst. Ser. B}, 14(1):41--74, 2010.

\bibitem{GonzalezHL13}
A.~Gonz\'alez, A.~Haro, and R.~de~la Llave.
\newblock Singularity theory for non-twist {KAM} tori.
\newblock {\em Mem. Amer. Math. Soc.}, 227(1067):vi+115, 2014.

\bibitem{GonzalezHL22}
A.~Gonz\'{a}lez, \`{A}. Haro, and R.~de~la Llave.
\newblock Efficient and reliable algorithms for the computation of non-twist
  invariant circles.
\newblock {\em Found. Comput. Math.}, 22(3):791--847, 2022.

\bibitem{HaroCFLM16}
A.~Haro, M.~Canadell, J.-Ll. Figueras, A.~Luque, and J.-M. Mondelo.
\newblock {\em The parameterization method for invariant manifolds}, volume 195
  of {\em Applied Mathematical Sciences}.
\newblock Springer, [Cham], 2016.
\newblock From rigorous results to effective computations.

\bibitem{HaroL19}
A.~Haro and A.~Luque.
\newblock A-posteriori {KAM} theory with optimal estimates for partially
  integrable systems.
\newblock {\em J. Differential Equations}, 266(2-3):1605--1674, 2019.

\bibitem{Henon66}
M.~H{\'e}non.
\newblock Exploration num\'erique du probl\`eme restreint iv. masses \'egales,
  orbites non p\'eriodiques.
\newblock {\em Bull. Astronom.}, 3(1--2):49--66, 1966.

\bibitem{Herman98}
M.-R. Herman.
\newblock { D\'emonstration d’un th\'eor\`eme de V.I. Arnold}.
\newblock S\'eminaire de Syst`emes Dynamiques et manuscript, 1998.

\bibitem{Kolmogorov54}
A.N. Kolmogorov.
\newblock On conservation of conditionally periodic motions for a small change
  in {H}amilton's function.
\newblock {\em Dokl. Akad. Nauk SSSR (N.S.)}, 98:527--530, 1954.
\newblock Translated in p. 51--56 of \emph{Stochastic Behavior in Classical and
  Quantum Hamiltonian Systems, Como 1977} (eds. G. Casati and J. Ford) Lect.
  Notes Phys. 93, Springer, Berlin, 1979.

\bibitem{Laskar05}
J.~Laskar.
\newblock Frequency map analysis and quasiperiodic decompositions.
\newblock In {\em Hamiltonian systems and {F}ourier analysis}, Adv. Astron.
  Astrophys., pages 99--133. Camb. Sci. Publ., Cambridge, 2005.

\bibitem{Laskar14}
J.~Laskar.
\newblock Michel h\'enon and the stability of the solar system.
\newblock In Hermann, editor, {\em Une vie d\'edi\'ee aux syst\`emes
  dynamiques: Hommage a Michel H\'enon}, pages 71--79. 2016.

\bibitem{LaskarR95}
J.~Laskar and P.~Robutel.
\newblock Stability of the planetary three-body problem. {I}. {E}xpansion of
  the planetary {H}amiltonian.
\newblock {\em Celestial Mech. Dynam. Astronom.}, 62(3):193--217, 1995.

\bibitem{LocatelliG05}
U.~Locatelli and A.~Giorgilli.
\newblock Construction of {K}olmogorov's normal form for a planetary system.
\newblock {\em Regul. Chaotic Dyn.}, 10(2):153--171, 2005.

\bibitem{LocatelliG07}
U.~Locatelli and A.~Giorgilli.
\newblock Invariant tori in the {S}un-{J}upiter-{S}aturn system.
\newblock {\em Discrete Contin. Dyn. Syst. Ser. B}, 7(2):377--398 (electronic),
  2007.

\bibitem{LuqueV09}
A.~Luque and J.~Villanueva.
\newblock Numerical computation of rotation numbers for quasi-periodic planar
  curves.
\newblock {\em Phys. D}, 238(20):2025--2044, 2009.

\bibitem{Moser62}
J.~Moser.
\newblock On invariant curves of area-preserving mappings of an annulus.
\newblock {\em Nachr. Akad. Wiss. G\"ottingen Math.-Phys. Kl. II}, 1962:1--20,
  1962.

\bibitem{Moser67}
J.~Moser.
\newblock Convergent series expansions for quasi-periodic motions.
\newblock {\em Math. Ann.}, 169:136--176, 1967.

\bibitem{Newton}
I.~S. Newton.
\newblock {\em Philosophiae naturalis principia mathematica}.
\newblock William Dawson \& Sons, Ltd., London, 1687.

\bibitem{Poincare_Book3BP}
H.~Poincar\'{e}.
\newblock {\em The three-body problem and the equations of dynamics}, volume
  443 of {\em Astrophysics and Space Science Library}.
\newblock Springer, Cham, 2017.
\newblock Poincar\'{e}'s foundational work on dynamical systems theory,
  Translated from the 1890 French original and with a preface by Bruce D. Popp.

\bibitem{PressTVF02}
W.H. Press, S.A. Teukolsky, W.T. Vetterling, and B.P. Flannery.
\newblock {\em Numerical Recipes in C: The Art of Scientific Computing}.
\newblock Cambridge University Press, second edition edition, 2002.

\bibitem{Robutel95}
P.~Robutel.
\newblock Stability of the planetary three-body problem. {II}. {KAM} theory and
  existence of quasiperiodic motions.
\newblock {\em Celestial Mech. Dynam. Astronom.}, 62(3):219--261, 1995.

\bibitem{Villanueva17}
J.~Villanueva.
\newblock A new approach to the parameterization method for {L}agrangian tori
  of {H}amiltonian systems.
\newblock {\em J. Nonlinear Sci.}, 27(2):495--530, 2017.

\end{thebibliography}

\appendix

\section{KAM Theorem and Iterative Lemma}\label{section: theorems}
Here we gather both the KAM Theorem and the Iterative Lemma in a tailored form. For a more
detailed exposition of them have a look at \cite{FiguerasHTheorem}.
\begin{theorem}\label{thm:KAM}
Let  $\H:\cmani\to \CC$ be a real-analytic Hamiltonian, defined in an open set $\cmani\subset \TT_\CC^m \times \CC^m$. 
Let $K:\bar\TT^m_\rho\to  \cmani$ be a continuous map, real-analytic in $\TT^m_\rho$,
whose derivatives are also continuous in $\bar\TT^m_\rho$, defining an homotopic to the zero-section  embedding of 
$\bar\TT^m_\rho$ into $\TT_\CC^m \times \CC^m$ (in particular $K(\theta)-(\theta,0)$ is $2\pi$-periodic).
 Let
$\omega\in \Dioph{\gamma}{\tau}{m}$ be a Diophantine vector, for some
$\gamma >0$ and $\tau \geq m-1$. 
We also assume:
\begin{itemize}[leftmargin=*]
\item [$H_1$]
There exist constants $\cteXH$, $\cteDXH$, $\cteDXHT$,
$\cteDDXH$ such that 
\[
\norm{\Xh}_{\cmani} \leq \cteXH, \quad
\norm{\Dif \Xh}_{\cmani} \leq \cteDXH,  \quad
\norm{(\Dif \Xh)^\ttop}_{\cmani} \leq \cteDXHT,  \quad
\norm{\Dif^2 \Xh}_{\cmani} \leq \cteDDXH.
\]
\item [$H_2$] 
There are  {\em condition numbers}  
$\sigmaDK$, $\sigmaDKT$, 
 $\sigmaB$, $\sigmaN$, $\sigmaNT$, and $\sigmaTinv$ such that
 \[
  \norm{\DK}_{\rho} < \sigmaDK, \quad
  \norm{\DKT}_{\rho} < \sigmaDKT,\quad
  \norm{\B}_\rho < \sigmaB,
  \]
  \[
\norm{\N}_{\rho} < \sigmaN, \quad
 \norm{\NT}_{\rho} < \sigmaNT, \quad
 \abs{\aver{T}^{\text{-}1}} < \sigmaTinv;
\]
\end{itemize}
Then, for each $\delta\in ]0,\rho/6[$,
there exists constants $\mathfrak{C}, \CtheoDeltaK$ depending on $\rho, \delta$ and the above constants and objects, such that, if 
\begin{equation}\label{eq:KAM:HYP}
\frac{\CtheoE}{\gamma\delta^{\tau+1}} \wnormeta < 1,
\end{equation}
where 
$
	\etaL=  - \NT\Omega\: (\Loper\K + \Xh\comp K),\ \etaN  =  \LT\Omega\: (\Loper\K + \Xh\comp K),
$
then, for $\rho_\infty= \rho-6\delta$, there exists $K_\infty:\bar\TT^m_{\rho_\infty}\to \cmani$
continuous, real-analytic in $\TT^m_{\rho_\infty}$,  
whose derivatives are also continuous in $\bar\TT^m_{\rho_\infty}$, defining an homotopic to the zero-section  embedding of 
$\bar\TT^m_{\rho_\infty}$ into $\cmani$ that is invariant under $\Xh$, with frequency $\omega$, so that
\[
\Loper\K_\infty + \Xh\comp \K_\infty= 0.
\]
Moreover, $\K_\infty$ satisfies hypothesis $H_2$, in $\TT^m_{\rho_\infty}$, and it is close to $K$:
\begin{equation*} 
\norm{\K_\infty-\K}_{\rho_\infty} \leq 
\frac{\CtheoDeltaK}{\gamma\delta^\tau} \wnormeta.
\end{equation*}
\end{theorem}

{The proof of the previous theorem consists of iteratively applying the following lemma.}

\begin{lemma}[The Iterative Lemma] \label{lem:iterative lemma}
Let us be under the same hypotheses as in Theorem \ref{thm:KAM}.
For any $\delta\in ]0,\rho/3]$, there exist constants
$\Csym$, 
$\CxiL$, $\CDeltaKn$, $\CDeltaDKn$, $\CDeltaDKnT$,
$\CDeltaBn$, $\CDeltaNn$, $\CDeltaNnT$, $\CDeltaiTn$, $\hCDelta$ and
$\QCetanL$, $\QCetanN$, 
such that if 
\begin{equation*}
\frac{\hCDelta}{\gamma \delta^{\tau+1}}  \wnormeta< 1,
\end{equation*}
then we have a new real-analytic parameterization $\Kn:\bar\TT^m_{\rho-2\delta}\to  \cmani$, 
that defines new objects $\Ln$, $\Bn$, $\Nn$ and $\Tn$ (obtained replacing $\K$ by $\Kn$ in the corresponding definitions) satisfying
\[
 \norm{\DKn}_{\rho-3\delta} < \sigmaDK,  
 \norm{\DKnT}_{\rho-3\delta} < \sigmaDKT,
 \norm{\Bn}_{\rho-3\delta} < \sigmaB,  
\]
\[
 \norm{\Nn}_{\rho-3\delta} < \sigmaN, 
 \norm{\NnT}_{\rho-3\delta} < \sigmaNT,  
 \abs{\aver{\Tn}^{\text{-}1}} < \sigmaTinv,  
\]
and
\[
 \norm{\Kn-\K}_{\rho-2 \delta} \leq 
\frac{\CDeltaKn}{\gamma\delta^\tau} \wnormeta,  
\]
\[
 \norm{\DKn-\DK}_{\rho-3 \delta} \leq 
\frac{\CDeltaDKn}{\gamma\delta^{\tau+1}} \wnormeta,  
\]
\[
 \norm{\DKnT-\DKT}_{\rho-3 \delta} \leq 
\frac{\CDeltaDKnT}{\gamma\delta^{\tau+1}} \wnormeta,  
\]
\[
 \norm{\Bn-\B}_{\rho-3\delta} \leq  
\frac{\CDeltaBn}{\gamma\delta^{\tau+1}} \wnormeta, 
\]
\[
 \norm{\Nn-\N}_{\rho-3 \delta} \leq 
\frac{\CDeltaNn}{\gamma\delta^{\tau+1}} \wnormeta,  
\]
\[
 \norm{\NnT-\NT}_{\rho-3 \delta} \leq 
\frac{\CDeltaNnT}{\gamma\delta^{\tau+1}} \wnormeta,  
\]
\[
 \abs{\aver{\Tn}^{\text{-}1}-\aver{\T}^{\text{-}1} } \leq  
\frac{\CDeltaiTn}{\gamma\delta^{\tau+1}} \wnormeta. 
\]
Moreover, the tangent and normal components of the new error of invariance 
\[
\En =  \Xh\comp \Kn + \Loper \Kn,
\]
satisfy
\begin{equation*}\label{lemma:CxinN}
\norm{\etanN}_{\rho-3\delta} \leq \frac{\QCetanN}{\delta} \qwnormeta,
\end{equation*}
and
\begin{equation*}\label{lemma:CxinL}
\begin{split}
\norm{\etanL}_{\rho-3\delta} 
& \leq \frac{\QCetanL}{\gamma\delta^{\tau+1}} \qwnormeta.
\end{split}
\end{equation*}
\end{lemma}

\end{document}